\documentclass[a4paper,12pt]{article}
\usepackage{amscd}
\usepackage{amsmath,amsfonts,amssymb,amscd}
\usepackage{indentfirst,graphicx,epsfig}
\usepackage{graphicx,psfrag}

\input{epsf}

\newtheorem{thm}{Theorem}

\newtheorem{problem}{Problem}
\newtheorem{lem}{Lemma}

\def\qed{\hfill \nopagebreak\rule{5pt}{8pt}}

\title{\bf Rainbow connection in $3$-connected graphs\footnote{Supported by NSFC and the
Fundamental Research Funds for the Central Universities. } }
\author{
\small  Xueliang Li, Yongtang Shi\\
\small Center for Combinatorics and LPMC-TJKLC \\
\small Nankai University, Tianjin 300071, China \\
\small lxl@nankai.edu.cn,  shi@nankai.edu.cn
\date{}}
\begin{document}
\maketitle
\begin{abstract}
An edge-colored graph $G$ is rainbow connected if any two vertices
are connected by a path whose edges have distinct colors. The
rainbow connection number of a connected graph $G$, denoted by
$rc(G)$, is the smallest number of colors that are needed in order
to make $G$ rainbow connected. In this paper, we proved that
$rc(G)\leq 3(n+1)/5$ for all $3$-connected graphs.
\\
[2mm] Keywords: rainbow connection; connectivity; the fan lemma\\
[2mm] AMS Subject Classification (2010): 05C15, 05C40.
\end{abstract}

\section{Introduction}
All graphs considered in this paper are simple, finite and
undirected. We follow the notation and terminology of Bondy and
Murty \cite{BM}. An edge-colored graph $G$ is {\it rainbow
connected} if any two vertices are connected by a path whose edges
have distinct colors. Obviously, if $G$ is rainbow connected, then
it is also connected. This concept of rainbow connection in graphs
was introduced by Chartrand et al. in \cite{CJMZ}. The {\it rainbow
connection number} of a connected graph $G$, denoted by $rc(G)$, is
the smallest number of colors that are needed in order to make $G$
rainbow connected. An easy observation is that if $G$ is of order
$n$ then $rc(G)\leq n-1$, since one may color the edges of one
spanning tree of $G$ with different colors and the remaining edges
with colors already used. It is easy to verify that $rc(G)=1$ if and
only if $G$ is a complete graph, that $rc(G)=n-1$ if and only if $G$
is a tree. Notice that for the cycle $C_n$ of order $n$,
$rc(C_n)=\lceil n/2\rceil$. It was shown that computing the rainbow
connection number of an arbitrary graph is NP-hard \cite{CLRTY}.

There are some approaches to study the bounds of $rc(G)$ with
respect to the minimum degree $\delta(G)$. In \cite{CLRTY} Caro et
al. have shown that if $G$ is a graph of order $n$ with minimum
degree $\delta$, then $rc(G)< \min \{(\ln
\delta/\delta)n(1+o_\delta(1)), (4\ln \delta+3)n/\delta\}$. By
employing the method of $2$-step dominating set, Krivelevich and
Yuster \cite{KY} have shown that a connected graph $G$ with $n$
vertices and minimum degree $\delta$ has $rc(G) < 20n/\delta$.
Schiermeyer \cite{S} proved that $rc(G) < 3n/ 4$ for graphs with
minimum degree three. Very recently, Chandran et al. \cite{CDRV}
have improved the upper bound of Krivelevich and Yuster by showing
that for every connected graph $G$ of order $n$ and minimum degree
$\delta$, $rc(G)\leq 3n/(\delta+1)+3$.

With respect to the the relation between $rc(G)$ and the
connectivity $\kappa(G)$, in \cite{S}, the author mentioned that
Hajo Broersma asked a question at the IWOCA workshop:
\begin{problem}
What happens with the value $rc(G)$ for graphs with higher
connectivity.
\end{problem}
Schiermeyer \cite{S} have shown that if $G$ is a graph of order $n$
with $\kappa(G)=1$ and $\delta\geq 3$, then $rc(G) \leq (3n-1)/ 4$.
In \cite{CLRTY} Caro et al. proved that if $\kappa(G)=2$ then $rc(G)
\leq 2n/ 3$. From the result of Chandran et al. \cite{CDRV}, we can
easily obtain an upper bound of the rainbow connection number:
$$rc(G)\leq \frac {3n} {\delta+1}+3\leq  \frac {3n} {\kappa(G)+1}+3.$$
Therefore, for $\kappa(G)=3$, $rc(G)\leq 3n/4+3$, and $\kappa(G)=4$,
$rc(G)\leq 3n/5+3$. In this paper, motivated by the results in
\cite{CLRTY}, we will improve this bound by showing the following
result.
\begin{thm}\label{thm}
If $G$ is a $3$-connected simple graph with $n$ vertices, then
$rc(G) \leq 3(n+1)/ 5$.
\end{thm}

Before proceeding, we recall the fan lemma, which will be used
frequently in the sequel.

\begin{lem}[The Fan Lemma] Let $G$ be a $k$-connected
graph, $x$ a vertex of $G$, and let $Y\subseteq V-\{x\}$ be a set of
at least $k$ vertices of $G$. Then there exists a $k$-fan in $G$
from $x$ to $Y$, namely there exists a family of $k$ internally
disjoint $(x, Y)$-paths whose terminal vertices are distinct in Y.
\end{lem}

\section{Proof of Theorem \ref{thm}.}

Let $H$ be a maximal connected subgraph of $G$ satisfying that
$rc(H)\leq 3h/5-1/ 5$, where $h$ is the number of vertices of $H$.

We first claim the existence of $H$. If $G$ contains a triangle
$C_3$, then we can choose the triangle as $H$, since $rc(C_3)=1<
8/5$. If $G$ contains $C_k$ ($k\geq 4$ and $k\neq 5$) as a subgraph,
then we take $H=C_k$, since $rc(C_k)=\lceil k/ 2\rceil\leq {3k}/
{5}-1/ 5$. Now suppose all the cycles contained in $G$ are of length
$5$, then we can take $H$ as the graph obtained by adding one
pendent edge to $C_5$. Observe that $h=6$ and $rc(H)=3<{17}/ {5}$.

We next claim that $h\geq n-3$. By contradiction. Suppose there are
four distinct vertices outside of $H$, denoted by $x_1,\, x_2,\,x_3,
\,x_4$. Then by the fan lemma, each of them has three internally
disjoint paths to $H$.

We assume first each of $x_1,\, x_2,\,x_3, \,x_4$ has three
neighbors in $H$. Let $f_{ij}$ be the edges incident to the vertex
$x_i$, $j=1,2,3$. We can add $x_1,\, x_2,\,x_3, \,x_4$ to $H$, and
form a lager subgraph $H'$ with $h+4$ vertices. Now we use only two
new colors $1$ and $2$ to color the $12$ edges. Assigning color $1$
to edges $f_{i1}$ for $i=1,2,3$ and color $2$ to other $9$ edges.
Then we have
$$rc(H')\leq rc(H)+2\leq {3h}/ {5}-1 /5 +2< 3(h+4)/5-1/ 5,$$
contradicting to the choice of $H$.

It follows that at least one of these four vertices, say $x$, has
the property that one of the three internally disjoint $(x,H)$-paths
$P_0,\,P_1,\,P_2$ has length at least two. Furthermore, among all
vertices satisfying the above property, we choose vertex $x$ such
that one of the three paths has length one, say $P_0=e_0$, and that
the sum of lengths of $P_1$ and $P_2$ is as large as possible.
Denote $P_1=au_1u_2\ldots u_sx$ and $P_2=xv_{1}v_{2}\ldots v_tb$
with $a,\,b\in H$ and $u_i,\,v_j\notin H$ for all $i$ and $j$. With
loss of generality, we assume $t\geq s$, and then $t\geq 1$. We
first assume $s+t\geq 3$. We can add
$v_1,\,v_2,\,\ldots,\,v_s,\,x,\,u_1,\,u_2,\,\ldots,\,u_t$ to $H$ and
form a larger subgraph $H'$ with $h+s+t+1$ vertices. If $s+t$ is
even, then we can color the $s+t+2$ edges of path $au_1u_2\ldots
u_sxv_{1}v_{2}\ldots v_tb$ with $(s+t+2)/2$ new colors. In the first
half of the path the colors are all distinct, and the same ordering
of colors is repeated in the second half of the path. We can color
edge $e_0$ with any color already appeared in $H$, and then it is
straightforward to verify that $H'$ is rainbow connected. If $s+t$
is odd, then we can color the $s+t+2$ edges of path $au_1u_2\ldots
u_sxv_{1}v_{2}\ldots v_tb$ with $(s+t+1)/2$ new colors as follows.
The middle edge of the path and edge $e_0$ get any color that
already used in $H$. The first $(s+t+1)/2$ edges of the path all
receive distinct new colors, and in the last $(s+t+1)/2$ edges of
the path this coloring is repeated in the same order. Again it is
straightforward to verify that $H'$ is rainbow connected. We now
have
\begin{align*}
rc(H')&\leq rc(H)+\left\lceil{(s+t+1)} /2\right\rceil\\
& \leq {3h}/{5} -1/ 5+\left\lceil {(s+t+1)}/ 2\right\rceil \leq
{3(h+s+t+1)}/{5}-1/ 5,
\end{align*}
contradicting the maximality of $H$. Hence, we only assume $1\leq
s+t\leq 2$. We consider three cases as follows.

\begin{figure}[ht]
\psfrag{x}{$x$}\psfrag{x1}{$x_1$}\psfrag{a}{$a$}
\psfrag{b}{$b$}\psfrag{v1}{$v_1$}\psfrag{v2}{$v_2$}\psfrag{u1}{$u_1$}
\psfrag{H}{$H$}\psfrag{1}{$1$}\psfrag{2}{$2$}\psfrag{3}{$3$}
\psfrag{a'}{$a'$}\psfrag{b'}{$b'$}\psfrag{v1'}{$v_1'$}\psfrag{v2'}{$v_2'$}\psfrag{u1'}{$u_1'$}
\begin{center}
\includegraphics[width=15cm]{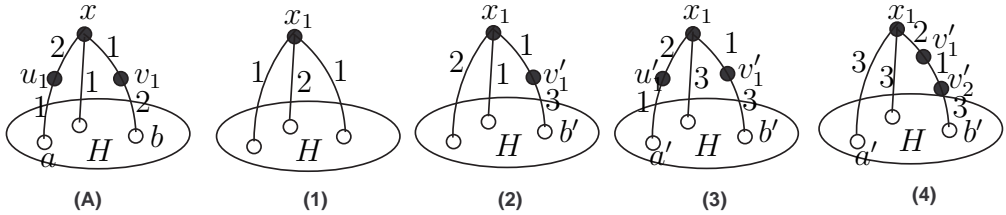}
\end{center}
\caption{$s+t=2$ and $P_1=au_1x$, $P_2=xv_1b$} \label{fig1}
\end{figure}

{\bf Case 1.} $s+t=2$ and $P_1=au_1x$, $P_2=xv_1b$ (see Figure
\ref{fig1}(A)).

\begin{figure}[ht]
\psfrag{x}{$x$}\psfrag{x1}{$x_1$}\psfrag{a}{$a$}
\psfrag{b}{$b$}\psfrag{v1}{$v_1$}\psfrag{v2}{$v_2$}\psfrag{u1}{$u_1$}
\psfrag{H}{$H$}\psfrag{c}{$c$}\psfrag{1}{$1$}\psfrag{2}{$2$}
\begin{center}
\includegraphics[width=8cm]{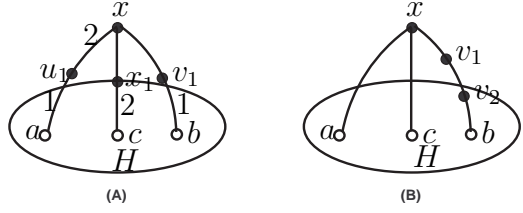}
\end{center}
\caption{Graphs used in Case 1 and Case 2.} \label{fig2}
\end{figure}

Since there are at least $4$ vertices outside of $H$, there exists
at least one vertex different from $x$, $u_1$ and $v_1$, say $x_1$.
By the choice of $x$, there is no $(x_1, x)$-path, $(x_1, u_1)$-path
and $(x_1, v_1)$-path without using any vertex of $H$ except one
case: there is one path of length two joining $x$ to $H$ through
$x_1$, say $P_3=xx_1c$ with $c\in H$. In this case, we only consider
the three paths $P_1$, $P_2$ and $P_3$  (as shown in Figure
\ref{fig2}(A)). We can add vertices $x$, $u_1$, $v_1$ and $x_1$ to
$H$ and form a larger subgraph $H'$ with $h+4$ vertices.
By assigning color $1$ to edges $au_1$, $bv_1$, color $2$ to edges
$u_1x$, $cx_1$, and one color already appeared in $H$ to edges
$v_1x$, $xx_1$, we have a contradiction as
$$rc(H')\leq rc(H)+2\leq {3h}/ {5}-1 /5 +2< 3(h+4)/5-1/ 5.$$
Now by the fan lemma, there are three internally disjoint
$(x_1,H)$-paths $P_0',\,P_1',\,P_2'$. By the choice of $x$, the
lengths of $P_0',\,P_1',\,P_2'$ only have four possibilities:

{\bf Subcase 1.1.} $1,\,1,\,1$. Let
$P_0'=e_0',\,P_1'=e_1',\,P_2'=e_2'$. We can add $x,\,u_1,\,v_1$ and
$x_1$ to $H$, and form a larger graph $H'$ of order $h+4$. By
assigning color $1$ to $au_1$, $e_0$, $xv_1$, $e_0'$, $e_1'$, and
color $2$ to $u_1x$, $v_2b$, $e_2'$, as shown in Figure
\ref{fig1}(1), we can obtain that $rc(H')\leq rc(H)+2\leq {3h}/
{5}-1/5 +2< 3(h+4)/5-1/ 5$, a contradiction.

{\bf Subcase 1.2.} $1,\,1,\,2$. Let
$P_0'=e_0',\,P_1'=e_1',\,P_2'=x_1v_1'b'$. We can add
$x,\,u_1,\,v_1$, $x_1$ and $v_1'$ to $H$, and form a larger graph
$H'$ of order $h+5$. By coloring all edges of paths
$P_0,\,P_1,\,P_2$ the same as Subcase 1.1 and assigning color $1$ to
$e_0'$, $x_1v_1'$, color $2$ to $e_1'$, and color $3$ to $v_1'b'$,
as shown in Figure \ref{fig1}(2), we can obtain that $rc(H')\leq
rc(H)+3\leq {3h}/ {5}-1/5 +3=3(h+5)/5-1/ 5$, a contradiction.

{\bf Subcase 1.3.} $1,\,2,\,2$. Let
$P_0'=e_0',\,P_1'=a'u_1'x_1,\,P_2'=x_1v_1'b'$. We can add
$x,\,u_1,\,v_1$, $x_1,\,u_1'$ and $v_1'$ to $H$, and form a larger
graph $H'$ of order $h+6$. By coloring all edges of paths
$P_0,\,P_1,\,P_2$ the same as Subcase 1.1 and assigning color $1$ to
$a'u_1'$, $x_1v_1'$, color $2$ to $u_1'x_1$, and color $3$ to
$e_0'$, $v_1'b'$, as shown in Figure \ref{fig1}(3), we can obtain
that $rc(H')\leq rc(H)+3\leq {3h}/ {5}-1/5 +3<3(h+6)/5-1/ 5$, a
contradiction.

{\bf Subcase 1.4.} $1,\,1,\,3$. Let
$P_0'=e_0',\,P_1'=e_1',\,P_2'=x_1v_1'v_2'b'$. We can add
$x,\,u_1,\,v_1$, $x_1,\,v_1'$ and $v_2'$ to $H$, and form a larger
graph $H'$ of order $h+6$. By coloring all edges of paths
$P_0,\,P_1,\,P_2$ the same as Subcase 1.1 and assigning color $1$ to
$v_1'v_2'$, color $2$ to $x_1v_1'$, and color $3$ to $e_0'$, $e_1'$,
$v_2'b'$, as shown in Figure \ref{fig1}(4), we can obtain that
$rc(H')\leq rc(H)+3\leq {3h}/ {5}-1/5 +3<3(h+6)/5-1/ 5$, a
contradiction.

{\bf Case 2.} $s+t=2$ and $P_1=ax$, $P_2=xv_1v_2b$ (see Figure
\ref{fig2}(B)).

Since $v_1\notin H$, there are three disjoint $(v_1,H)$-paths by the
fan lemma. Then there is at least one additional $(v_1,H)$-path
$P_3$ except paths $v_1xa$ and $v_1v_2b$. By the choice of $x$, the
length of $P_3$ must be at most two. If $P_3$ is of length two, then
this is the case of Figure \ref{fig2}(A), we have done. If $P_3$ is
of length one, then paths $axv_1$, $v_1v_2b$ and $P_3$ build the
same structure as Case 1, and thus we have done.

{\bf Case 3} $s+t=1$.

Since $t\geq s$, we have $t=1$. Now we can assume $P_1=e_1$ and
$P_2=xv_1b$. Then there are at least two distinct vertices outside
of $H$ different from $x$ and $v_1$, say $x_1$ and $x_2$. Similarly,
for $i=1,\,2$, there is no $(x_i, x)$-path and $(x_i, v_1)$-path
without using any vertex of $H$. So there are also three internally
disjoint $(x_1,H)$-paths $P_0',\,P_1',\,P_2'$ and $(x_2,H)$-paths
$P_0'',\,P_1'',\,P_2''$, respectively. If all these paths are of
length one, then we can add $x,\,v_1$, $x_1$ and $x_2$ to $H$, and
form a larger graph $H'$ of order $h+4$. By assigning color $1$ to
edges $e_0$, $xv_1$, $P_0',\,P_1',\,P_0'',\,P_1''$, color $2$ to
edges $e_1$, $v_1b$, $P_2',\,P_2''$, we can obtain that $rc(H')\leq
rc(H)+2\leq {3h}/ {5}-1/5 +2<3(h+4)/5-1/ 5$, a contradiction.
Otherwise, without loss of generality, we assume one of the three
$(x_1,H)$-paths $P_0',\,P_1',\,P_2'$ has length $2$. Let
$P_0'=e_0',\,P_1'=e_1',\,P_2'=x_1v_1'b'$. We can add $x,\,v_1$,
$x_1$ and $v_1'$ to $H$, and form a larger graph $H'$ of order
$h+4$. By assigning color $1$ to edges $e_0$, $e_1$, $xv_1$, $v_1b$,
color $2$ to edges $e_0'$, $e_1'$, $x_1v_1'$, $v_1'b'$, we can
obtain that $rc(H')\leq rc(H)+2\leq {3h}/ {5}-1/5 +2<3(h+4)/5-1/ 5$,
a contradiction.

Now we have proved that $h\geq n-3$. By considering some cases, we
can easily obtain that $rc(G)\leq 3(n+1)/5$: if $h=n-3$, then
$rc(G)\leq rc(H)+2\leq 3(h-3)/5-1/5+2<3(n+1)/5$; if $h=n-2$, then
$rc(G)\leq rc(H)+2\leq 3(h-2)/5-1/5+2=3(n+1)/5$; if $h=n-1$, then
$rc(G)\leq rc(H)+1\leq 3(h-1)/5-1/5+1<3(n+1)/5$.

The proof is completed.\qed

\end{document}